\documentclass[a4paper,10pt]{article}
\usepackage{amsfonts, amsmath, cite}

\newcommand{\cF}{\mathcal{F}}
\newcommand{\R}{\mathbb{R}}
\newcommand{\Z}{\mathbb{Z}}
\newcommand{\N}{\mathbb{N}}
\renewcommand{\H}{\mathbb{H}}
\newcommand{\cM}{\mathcal{M}}
\newcommand{\cU}{\mathcal{U}}
\newcommand{\cV}{\mathcal{V}}
\newcommand{\noi}{\noindent}
\newtheorem{theorem}{Theorem}
\newtheorem{lemma}{Lemma}
\newtheorem{proposition}{Proposition}
\newtheorem{corollary}{Corollary}

\newtheorem{example}{Example}
\newtheorem{remark}{Remark}

\newenvironment{proof}{\smallskip \noindent{\bf Proof: }}{\hfill$\Box$\medskip}

\title{Horocycle flows for laminations by hyperbolic Riemann surfaces and Hedlund's theorem}

\author{Matilde Mart\'{\i}nez\thanks{Instituto de Matem\'atica y Estad\'{\i}stica Rafael Laguardia, Facultad de
Ingenier\'{\i}a,
Universidad de la Rep\'ublica. J.Herrera y Reissig 565,C.P.11300, Montevideo, Uruguay. matildem@fing.edu.uy.
Research partially supported by ANII FCE-2007-106 and
CSIC , Proyecto de Grupo de Investigación 618 (Sistemas Din\'amicos} \and Shigenori Matsumoto
\thanks{ Department of Mathematics, College of
Science and Technology, Nihon University, 1-8-14 Kanda, Surugadai,
Chiyoda-ku, Tokyo, 101-8308 Japan, matsumo@math.cst.nihon-u.ac.jp.
Research partially supported by Grant-in-Aid for
Scientific Research (C) No.\ 25400096.}
\and
Alberto Verjovsky\thanks{Universidad Nacional Aut\'onoma de M\'exico.
Apartado Postal 273. Admon. de correos \#3, C.P. 62251 Cuernavaca,
Morelos, Mexico. alberto@matcuer.unam.mx. Research partially
supported by CONACyT proyecto U1 55084, and PAPIIT (Universidad
Nacional Aut\'onoma de M\'exico) \# IN102108.} }

\begin{document}

\maketitle

\begin{center}
{\it To \'Etienne Ghys, on the occasion of his 60th birthday.}
\end{center}

\begin{abstract}
We study the dynamics of the geodesic and horocycle flows of the unit tangent bundle $(\hat M, T^1\cF)$ of a compact
minimal
lamination $(M,\mathcal F)$ by negatively curved surfaces. We give conditions under which the action of the affine group
generated by the joint action of these flows is minimal, and examples where this action is not minimal.
In the first case, we prove that if $\cF$ has a leaf which is not simply connected,
the horocyle flow is topologically transitive.
\end{abstract}

\section{Introduction}

The geodesic and horocycle flows over compact hyperbolic surfaces
have been studied in great detail since the pioneering work in the
1930's by E. Hopf and G. Hedlund. Such flows are particular
instances of flows on homogeneous spaces induced by one-parameter
subgroups, namely, if $G$ is a Lie group, $K$ a closed subgroup and
$N$ a one-parameter subgroup of $G$, then $N$ acts on the
homogeneous space $K\backslash{G}$ by right multiplication on left
cosets. One very important case is when $G=SL(n,\R)$,
$K=SL(n,\Z)$ and $N$ is a unipotent one parameter subgroup of
$SL(n,\R)$, i.e., all elements of $N$ consist  of matrices
having all eigenvalues equal to one. In this case
$SL(n,\Z)\backslash{SL(n,\R)}$ is the space of unimodular lattices. By
a theorem by Marina Ratner (see \cite{Ratner}), which gives a positive answer to the
Raghunathan conjecture, the closure of the orbit under the unipotent
flow of a point $x\in{SL(n,\Z)\backslash{SL(n,\R)}}$ is the
orbit of $x$ under the action of a closed subgroup $H(x)$. This
particular case already has very important applications to number
theory, for instance, it was used by G. Margulis and Dani in \cite{Dani-Margulis} and
Margulis in \cite{Margulis} to give a positive answer to the Oppenheim conjecture. When
$n=2$ and $\Gamma$ is a discrete subgroup of $SL(2,\R)$ such
that $M:=\Gamma\backslash{SL(2,\R)}$ is of finite Haar volume,
and $N$ is any unipotent one-parameter subgroup acting on $M$,
Hedlund proved that any orbit of the flow is either a
periodic orbit or dense. When $\Gamma$ is cocompact the flow induced
by $N$ has every orbit dense, so it is a minimal flow. The
horocycle flow on a compact hyperbolic surface is a homogeneous
flow of the previous type and most of the important dynamic,
geometric and ergodic features of the general case are already present in this
3-dimensional case.

On the other hand, the study of Riemann surface laminations has
recently played an important role in holomorphic dynamics (see \cite{Ghys2} and \cite{Lyubich-Minsky}), polygonal
tilings of the Euclidean or hyperbolic plane (see \cite{Bellissard-Benedetti-Gambaudo}, \cite{Petite}), moduli spaces of
Riemann surfaces (see \cite{Nag-Sullivan}), polygonal billiards (see \cite{McMullen} and \cite{Zorich}), etc. It is
natural then to consider compact
laminations by surfaces with a Riemannian metric of negative
curvature and consider the positive and negative horocycle flows on
the unit tangent bundle of the lamination. In the spirit of Raghunathan and Ratner, in this paper we study the closures
of the horocycle orbits in this
non-transitive $PSL(2,\R)$-space. We give examples of orbit closures which are not algebraic in the sense of Ratner.
We give sufficient conditions for the horocycle flow
to be minimal. Finally, we state a conjecture relating the minimality of the joint action of the geodesic and horocycle
flows and the minimality of the horocycle flow alone.

The authors would like to thank \'E.Ghys for explaining the argument that became the proof of Proposition
\ref{proposition:more_than_one_point_implies_minimality}
and B.Deroin for valuable insights during informal discussions of this work. We would also like to thank the referees
at {\it Journal of Modern Dynamics} for a careful reading of the paper and useful suggestions that contributed to its
improvement.


\section{Laminations by hyperbolic surfaces}

For the definition of a compact {\em lamination} or {\em foliated
space} in general, we refer the reader to Chapter 11 of \cite{Candel-Conlon}. In this paper, we only concern ourselves
 with compact laminations by
hyperbolic surfaces.

\subsection{Laminations by hyperbolic surfaces and their unit tangent bundles}
A compact {\em lamination by hyperbolic surfaces} $(M,\cF)$ (or simply $\cF$) consists of a compact metrizable space $M$
together
with a family $\{(U_\alpha, \varphi_\alpha)\}$
such that
\begin{itemize}
\item[\textbullet] $\{U_\alpha\}$ is an open covering of $M$,
\item[\textbullet] $\varphi_\alpha: U_\alpha \to D_\alpha\times T_\alpha$ is a
	     homeomorphism, where $D_\alpha$ is a bounded open disk in
the Pioncar\'e upper half plane $\H$ and $T_\alpha$
is a topological space, and
\item[\textbullet] for $(y,t)\in \varphi_\beta(U_\alpha\cap U_\beta)$,
$\varphi_\alpha \circ \varphi_\beta^{-1}(y,t) = (g_{\alpha\beta}(t)y, h_{\alpha\beta}(t))$, where
$g_{\alpha\beta}(t)\in PSL(2,\R)$.
\end{itemize}
Notice that since $\varphi_\alpha$
is a homeomorphism, the map $h_{\alpha\beta}$ is a homeomorphism
from an open subset of $T_\beta$ to an open subset of $T_\alpha$,
and $g_{\alpha\beta}(t)$
is continuous in $t$.
In the sequel $M$ will always be compact.
Each $U_\alpha$ is called a {\em foliated chart}, a set of the form $\varphi_\alpha^{-1}(\{y\}\times T_\alpha)$ being its
{\em transversal}.
The sets of the form $\varphi_\alpha^{-1}(D_\alpha\times\{t\})$, called {\em plaques}, glue together to form maximal
2-dimensional connected manifolds called {\em leaves}. The lamination
$\cF$ is the partition of $M$ into leaves.
The leaves have structures of hyperbolic surfaces with local
charts the restriction of $p_1\circ\varphi_\alpha$'s to plaques,
where $p_1:D_\alpha\times T_\alpha\to D_\alpha$
is the projection onto the first factor.

A lamination is said to be {\em minimal} if all the leaves are dense.

A {\em holonomy invariant} measure is a family of Radon measures in all transversals which are invariant
under the maps $h_{\alpha\beta}$ when restricted to their domains and
ranges.

\bigbreak
For each $x\in M$, let us denote by $T^1_x\cF$ the set of the unit
tangent vectors of the leaf through $x$, and define the {\em laminated unit
tangent bundle} $\hat M$ of $(M,\cF)$ by
 $\hat M=\cup_{x\in M}T_x^1\cF$.
Points of $\hat M$ are denoted by $(x,v)$, where $x\in M$ and $v\in T^1_x\cF$.
Denote the natural projection by
$\Pi:\hat M\to  M$. Thus $\Pi^{-1}(x)=T^1_x\cF$.
There are homeomorphisms
$\psi_\alpha:\Pi^{-1}(U_\alpha)\to T^1D_\alpha\times T_\alpha$
defined by the leafwise differentials of
$\phi_\alpha$.
The locally defined transition function
$\psi_\alpha\circ\psi_\beta^{-1}:T^1D_\beta\times T_\beta\to
T^1D_\alpha\times T_\alpha$
is given by
$$\psi_\alpha \circ \psi_\beta^{-1}(w,t) = (g_{\alpha\beta}(t)_*(w),
h_{\alpha\beta}(t)),$$
where $g_{\alpha\beta}(t)_*$ is the differential of $g_{\alpha\beta}(t)$.
Plaques $\psi_\alpha^{-1}(T^1D_\alpha\times\{t\})$ all together
define a lamination which we shall denote by $T^1\cF$.

When we identify $T^1\H$ with $PSL(2,\R)$ as usual, then
$T^1D_\alpha$ and $T^1D_\beta$ are open subsets of $PSL(2,\R)$, and
the transition function takes the form
$$\psi_\alpha \circ \psi_\beta^{-1}(g,t) = (g_{\alpha\beta}(t) g, \,
h_{\alpha\beta}(t)),$$
where $g\in PSL(2,\R)$ and
$g_{\alpha\beta}(t) g$ is the group multiplication.
Since the
transition functions commute with the locally defined right translation of
$PSL(2,\R)$ on the first factor, we obtain a right action of the universal covering group
of $PSL(2,\R)$ on $\hat M$. Since the rotation by $2\pi$ is the identity on each
foliated chart, we
actually obtain a right action of $PSL(2,\R)$. The orbit foliation
of this action is just $T^1\cF$.

\subsection{Laminated geodesic and horocycle flows}
We are particularly interested in the right action of the following three
subgroups $D$, $U$
and $B$ of $PSL(2,\R)$, where

$$D=\left\{
\left (
\begin{array}{cc}
e^\frac{t}{2} & 0\\
0 & e^{-\frac{t}{2}}
\end{array}
\right );\ t\in\R
\right \},\ \
U=\left\{
\left (
\begin{array}{cc}
1 & t\\
0 & 1
\end{array}
\right );\ t\in\R
\right \}
\mbox{ and }
$$

$$B=\left\{
\left (
\begin{array}{cc}
a & b\\
0 & a^{-1}
\end{array}
\right );\ a>0,\ b\in\R
\right \}.$$

The right action of $D$ (resp.\ $U$) is called the laminated {\em geodesic}
(resp.\ {\em horocycle}) flow and denoted by $g^t$ (resp.\ $h^t$).
They satisfy
$$
g^t\circ h^s\circ g^{-t}=h^{se^{-t}}.
$$
For each leaf $L$ of $\cF$, the unit tangent bundle $T^1L$
is a leaf of $T^1\cF$, and the flows $g^t$ and $h^t$ restricted
to $T^1L$ are
the usual geodesic and horocycle flows of the hyperbolic surface $L$.


To study further the right $B$-action,
we consider the following
 identification  of $T^1\H$ with $\H\times{\mathbb S}^1$,
where $\mathbb S^1$ is the circle at infinity of $\H$.
For any
$w\in T^1\H$, denote by $w(\infty)$ the point in $\mathbb S^1$ at which
the geodesic with initial vector $w$ terminates, and by
 $p:T^1\H\to \H$ the bundle projection.
Now we identify $w\in T^1\H$ with
$(p(w),w(\infty))\in \H\times\mathbb S^1$.
For any  $g\in PSL(2,\R)$, $w\in T^1\H$,
we have $(p(g_*(w)),(g_*(w))(\infty))=(gp(w),gw(\infty))$,
where $gw(\infty)$ is the action of $PSL(2,\R)$ on $\mathbb S^1$.

Thus when $T^1D_\alpha\times T_\alpha$
and $T^1D_\beta\times T_\beta$ are identified with
$D_\alpha\times\mathbb S^1\times T_\alpha$
and  $D_\beta\times\mathbb S^1\times T_\beta$,
the transition function
$\psi_\alpha\circ\psi_\beta^{-1}$ takes the form
$$
\psi_{\alpha}\circ\psi_\beta^{-1}(y,\theta,t)=(g_{\alpha\beta}(t)y,g_{\alpha\beta}(t)\theta,h_{\alpha\beta}(t)).$$
On the other hand, we have $wb(\infty)=w(\infty)$ for any
$w\in T^1\H$ and $b\in B$.
Therefore the local right $B$ translation on a lamination chart
$D_\alpha\times \mathbb S^1\times T_\alpha$ leaves the second and the
third factor invariant.
In what follows, we shall depress the map $\phi_\alpha$ and
the subscript $\alpha$, and write $D\times\mathbb S^1\times T$
to denote a laminated chart in $\hat M$. The foliation of $\hat M$
obtained by the plaques
$D\times\{\theta\}\times\{t\}$ is the orbit foliation of the
right $B$-action.


\subsection{Abundance of compact laminations by hyperbolic surfaces}

Let $(M,\cF)$ be a compact 2-dimensional lamination
equipped with a leafwise Riemannian metric which is continuous and
leafwise $C^2$.
The restriction to a leaf $L$ of the metric lifts to a Riemannian metric on the universal cover $\tilde{L}$ of $L$.
Let us fix $x\in \tilde{L}$ and denote $A(x,r)$ the area of the ball on
$\tilde{L}$ centered at $x$ with radius $r>0$.
Consider the following condition (\ref{e-m3}).
\begin{equation}\label{e-m3}
  \lim_{r\to\infty} \frac{1}{r}\log A(x,r) >0,\ \ \forall x\in M.
\end{equation}

This condition is independent of the Riemannian metric chosen,
since all the leafwise metrics on the compact lamination $(M,\cF)$
are quasi-isometric on leaves.
A theorem due to Candel and Verjovsky (see \cite{Candel} and
\cite{Verjovsky}) guarantees
the existence of a Riemannian metric
of constant curvature $-1$ in the same conformal class of the metric we
started with. Namely we have the following theorem.

\begin{theorem}
 If $(M,\cF)$ is a compact 2-dimensional lamination which satisfies {\em (\ref{e-m3})}, then $\cF$ has a structure of a lamination by hyperbolic
surfaces.
\end{theorem}

\bigbreak

Such laminations which satisfy (\ref{e-m3}) are quite abundant. For
example, laminations without holonomy invariant measures are known to
satisfy (\ref{e-m3}).


\subsection{Harmonic measures}

The study of the laminated Brownian motion was developed by Lucy Garnett
in \cite{Garnett}. If $(M,\cF)$ is a compact lamination by hyperbolic
surfaces, the laminated Brownian motion
starting at $x\in M$ is defined simply to be the usual
Brownian motion on the leaf $L$ of $\cF$ through $x$
with respect to the hyperbolic metric on $L$.
A probability measure of $M$ which is stationary for this process is called
{\em harmonic}. A harmonic measure $\mu$ is {\em ergodic}
if every measurable subset of $M$ which is a union of leaves of $\cF$ has either $\mu$ measure zero or one.

Lucy Garnett proved the following ergodic theorem for harmonic measures:

\begin{theorem}[Lucy Garnett]
Let $\mu$ be an ergodic harmonic measure on $(M,\cF)$ and $f\in L^1(\mu)$. Then for $\mu$-almost every $x\in M$ and almost every
Brownian path $\omega(t)$ on the leaf through $x$
starting at $x$
\begin{equation}
\lim_{T\to\infty}\frac{1}{T}\int_0^T f(\omega(t))\, dt =\int_M f\, d\mu.
\end{equation}
\end{theorem}

When the lamination $(M,\cF)$ has a holonomy invariant measure $\nu$, we can consider the measure on $M$ which is
locally the product of $\nu$ and the volume on leaves. This gives a harmonic measure, and harmonic measures of this
kind are called {\em completely invariant}.


\section{The central stable foliation of the laminated geodesic flow}

As before, let $(M,\cF)$ be a compact minimal lamination by hyperbolic surfaces. The laminated geodesic
flow on $(\hat M, T^1\cF)$ is
an Anosov flow on
each leaf $T^1L$,
whose central stable foliation is the foliation by orbits of the action of the affine group $B$ restricted to $T^1L$.
Without any
further assumption on the dynamics of
the foliation, the central-stable manifolds of the laminated geodesic flow might be bigger than these $B$-orbits.
Nevertheless, we will refer
to the foliation by $B$-orbits on $\hat M$
as the {\em central-stable foliation of the laminated geodesic flow}, and to its leaves as the {\em central-stable leaves}.

Let $\Pi:\hat M \to M$ be the canonical projection $\Pi(x,v)=x$; for every $x\in M$ $\Pi^{-1}(\{x\})=T^1_x\cF$.
Under this projection, each central-stable leaf in $\hat M$ projects onto a leaf of $\cF$.

Let $\cM\subset\hat M$ be a compact minimal set for the central-stable foliation. Since $\cF$ is minimal and
central-stable leaves project
onto entire $\cF$-leaves, $\Pi(\cM)=M$.
Namely, $\cM\cap T^1_x\cF$ is never empty for $x\in M$. We will see that the number of points in this intersection
carries crucial
information on
the dynamics of the $B$-action.

It is possible that $\cM\cap T^1_x\cF$ consists of a single point for all $x$. Foliations for which this happens are characterized in the following Proposition.

\begin{proposition}
\label{proposition:one_point_implies_affine_action}
Let $(M,\cF)$ be a compact minimal lamination by hyperbolic surfaces.
Then $(M,\cF)$ is the lamination by orbits of a continuous, locally free, right action of the group
$B$ on $M$ if and only if there is a compact set
$\mathcal{M}\subset\hat M$
which is minimal for the action of $B$ on $\hat M$ such that $T^1_x\cF \cap \mathcal{M}$ consists of a single point
for all $x\in M$.

\end{proposition}

\begin{proof}
Assume there exists one such $\mathcal{M}$. The fact that it intersects each unit tangent space in only one point means that
$\pi:\mathcal{M}\to M$ is a homeomorphism. Since it
sends central-stable leaves in $\mathcal{M}$ to leaves of $\cF$, via $\Pi$ the right $B$-action on $\mathcal{M}$  defines a
right $B$-action on $M$ whose orbits are the leaves of $\cF$.

To show the converse, we need to study the geometry associated to a
locally free $B$-action. First of all, we fix a left invariant
hyperbolic Riemannian metric $g_B$ on $B$. In fact, the group $B$ is
isomorphic to
$$
B'=\left\{\left(\begin{array}{cc}
y&x\\
0&1
\end{array}\right)\mid x\in\R,y>0\right\},$$
and the metric
$$
g_{B'}=y^{-2}(dx^2+dy^2)$$
is left invariant. We choose $g_B$ to be the metric on $B$ which
corresponds to $g_{B'}$ by the isomorphism. Let $\phi:M\times B\to M$ be
a locally free $B$-action. For any $x\in M$, define $\phi_x:B\to M$ by
$$\phi_x(b)=\phi(x,b).$$
Then there is a leafwise metric $g_\phi$ of the orbit foliation which
satisfies
$$
\phi_x^*g_\phi=g_B,\ \mbox{ for any }\ x\in M.$$
The metric $g_\phi$ is well defined since $g_B$ is left invariant.

Now assume that the lamination $\cF$ on $M$ is defined by a right $B$
action $\phi$. Let $D$ and $U$ be the one-parameter subgroups of $B$
defined in the previous section. On the universal covering space of
each $B$-orbit (equipped with the metric $g_\phi$) of a point
$x\in M$, there is a point
$\xi(x)$ at infinity such that any $D$-orbit is a geodesic tending to $\xi(x)$
and any $U$-orbit is a horocycle at $\xi(x)$. Now since $M$ is compact,
the leafwise metric $g_\phi$ is quasi-isometric to the
leafwise metric $g$ which is given a priori. Therefore the point
$\xi(x)$ at infinity with respect to $g_\phi$ corresponds to a point
at infinity with respect to $g$, which we still denote by $\xi(x)$.
Notice that $\xi(x)$ is $B$-invariant and continuous in $x$.

Finally for any $x\in M$, define $X(x)\in T_x^1{\cF}$
to be the tangent vector tending to $\xi(x)$.
This way, we get a cross section $X:M\to\hat M$. The image $\mathcal M$
is compact, invariant and minimal for the action of $B$ on $\hat M$

%

\end{proof}

Foliations that arise
in this way have been widely studied, see for example \cite{Plante}. They do not possess holonomy invariant measures.
All of their leaves are homeomorphic to planes or cylinders, and in codimension one they must have countably many cylindrical leaves,
and these leaves have non-trivial holonomy.

Although the $C^\infty$ locally free $B$-actions in dimension 3 are found only in the $S^1$-bundles over surfaces (\cite{Asaoka}),
continuous actions are abundant. The central stable foliation of any topologically transitive Anosov flow is the orbit foliation of a continuous $B$-action, thanks to the existence of the Margulis measure along the stable foliation (\cite{Katok-Hasselblatt}). There are
many examples of such flows (\cite{Fenley}, \cite{Barbot}, \cite{Foulon-Hasselblatt}).

\bigbreak

Apart from the situation described in Proposition \ref{proposition:one_point_implies_affine_action}, two possibilities remain:
that the intersection $T^1_x\cF \cap \mathcal{M}$ has more than one point for all $x$, and that
$T^1_x\cF \cap \mathcal{M}$ has one point for some values of $x$ and more than one point for another values of $x$.

\begin{proposition}
\label{proposition:more_than_one_point_implies_minimality}
Let $(M,\cF)$ be a compact minimal lamination by hyperbolic surfaces.
If $\mathcal{M}\subset \hat M$ is compact and minimal for the action of $B$ on $\hat M$ and it intersects each unit tangent
space in more than one point, then $\mathcal{M} = \hat M$.
\end{proposition}

To prove this we will use an argument due to \'Etienne Ghys.

\begin{proof}
Assume $\mathcal{M}$ is a nontrivial minimal set (that is, that $\mathcal{M}$ is not the whole $\hat M$),
and that it intersects each unit tangent space in at least two points.

Let us take a point $x\in M$ and call $L_x$ the leaf of $\cF$ that passes through $x$. The fact that the lamination $(\hat M, T^1\cF)$ is minimal and $\mathcal{M}$
is $B$-invariant implies that
the unit tangent space to $L_x$ at $x$ is not contained in $\mathcal{M}$, i.e. $T^1_x\cF\cap \mathcal{M}$ is a nonempty closed proper subset of $T^1_x\cF$.
In any case, by identifying in the usual way $T^1_x\cF$ with the set of points at infinity of the universal cover $\tilde L_x$ of $L_x$,
we may think of $T^1_x\cF\cap \mathcal{M}$ as a subset $K_x$
of the circle at infinity of $\tilde L_x$. Notice that $K_x$ does not depend on $x$ but only on the leaf $L_x$, since $\mathcal{M}$ is $B$-invariant.

For the moment we will also assume that all leaves of $\cF$ are simply connected.

For every $x\in M$, let $\hat{K_x}$ be the convex hull of $K_x$ in $L_x$. (Namely, consider all geodesics in $L_x$ joining pairs of points in $K_x$,
and take the convex hull of their union. This set is $\hat{K_x}$.) It is possible to do this because $K_x$ has at least two points.

Let $f:M\to [0,+\infty)$ be the function defined by
$$
f(y)=d(y,\hat{K_y}),
$$
where $d$ is the hyperbolic distance on the leaf passing through $y$.

\begin{remark}
The function $f$ is measurable.
\end{remark}

Proof of the remark:  Let $E\simeq D\times T$ be a compact foliated
chart of $(M,\cF)$, where $D$ is a closed disk in $\R^2$ and $T$ is a
topological space. In $E$,  the $\hat{K_x}$  form a semicontinuous
family of compact sets parametrized by $T$, and the function $f$ is
the distance on each fiber $D\times \{t\}$ to the corresponding
compact set. Therefore $f$ is measurable.~$\bigtriangleup$

In \cite{Buber-Kirk}, the authors prove that the axiom of choice is not needed to prove the existence of minimal sets.
This implies that the function $f$ can be defined without using the axiom of choice, and this is another reason why we
can safely assume that $f$
is measurable (see \cite{Solovay}).

\bigbreak

Let $\mu$ be an ergodic harmonic measure on $(M,\cF)$. For every $n\in\N$, we define
$$
A_n=\{x\in M:\ f(x)\leq n\}.
$$
The sequence $\{A_n\}$ is increasing and $\mu(\cup_n A_n)=1$, therefore there exists an $n\in\N$ for which $\mu(A_n)>0$.
The ergodic theorem tells us that for
$\mu$-almost every $x\in M$ and almost every continuous path $\omega$ on $L_x$ that starts at $x$
\begin{equation}
\label{eq:behaviour_of_almost_every_path}
\lim_{T\to\infty} \frac{1}{T} \int_0^T \chi_{A_n}(\omega(t))\, dt>0,
\end{equation}
where $\chi_{A_n}$ is the characteristic function of the set $A_n$.

In spite of this, for every $x\in M$, the set of continuous paths $\omega(t)$ on $L_x$ which start at $x$ and which converge,
when $t\to\infty$, to a point outside $K_x$ has positive Wiener measure, and for any of these paths
$$
\lim_{T\to\infty}\frac{1}{T}\int_0^T \chi_{A_n}(\omega(t))\, dt =0.
$$
Since this contradicts (\ref{eq:behaviour_of_almost_every_path}), we have proved the lemma when all leaves of $\cF$ are hyperbolic planes.

If there is a leaf $L_x$ of $\cF$ which is not simply connected, we can still define $\hat{K_x}$ on its universal cover.
For points on $L_x$, we define $f(y)$ as the distance from $y$ to the projection of $\hat{K_x}$, and use
the same argument as before to complete the proof.
\end{proof}

In fact, the previous argument proves the slightly stronger statement, which will be useful later on:

\begin{proposition}
\label{proposition:generic_property_for_harmonic_measures}
Let $(M,\cF)$ be a compact minimal lamination by hyperbolic surfaces, and let $\mu$ be an ergodic harmonic measure on $M$.
If $\mathcal{M}\subset \hat M$ is compact and minimal for the action of $B$ on $\hat M$, then either $T^1_x\cF \cap \mathcal{M}$ is
a singleton for $\mu$-almost every $x$ or $\mathcal{M}=\hat M$.
\end{proposition}

\bigbreak

\begin{remark}
\label{remark:minimality_not_required}
Notice that in the proof of Proposition \ref{proposition:more_than_one_point_implies_minimality}, we do not actually use that $\mathcal M$
is minimal, only that it is a closed $B$-invariant set.
\end{remark}

We use this fact to prove the uniqueness of $B$-minimal sets, even in the case where the $B$ action is not minimal.

\begin{theorem}
\label{theorem:uniqueness_of_B_minimal_set}
 If the action of $B$ on the unit tangent bundle $\hat{M}$ is not minimal, then there exists a unique minimal set
$\cM$ which is a proper subset of $\hat{M}$.
\end{theorem}

\begin{proof}
 If $\cM_1$ and $\cM_2$ were two minimal sets they have to be disjoint and then $\cM_1\cup\cM_2$ would be a $B$-invariant set
which intersects the circle fibers of $\hat{M}$ in more than one point and therefore, by Proposition 2, $\cM_1\cup\cM_2=\hat{M}$ which is a contradiction.
\end{proof}

Proposition \ref{proposition:more_than_one_point_implies_minimality}, in the slightly stronger version mentioned in Remark \ref{remark:minimality_not_required}, has a pair of interesting consequences for the action of the affine group $B$, as we will now see.

\begin{theorem}
\label{theorem:transitivity_of_B}
If $M$ is a compact {\emph minimal} lamination by hyperbolic surfaces and $\hat{M}$ is the laminated
unit tangent bundle then the locally free action of $B$ on  $\hat{M}$ is topologically transitive. In particular the set of $B$-orbits
which are dense in  $\hat{M}$ is residual (in the sense of Baire).
\end{theorem}

\begin{lemma} Let $\mathcal{U}\subset\hat{M}$ be any open set. Let $B\cU=\{b(u) |\,b\in B, u\in \cU\}$ be the $B$-orbit of $\cU$. Then
$\overline{B\cU}=\hat{M}$, in other words $B\cU$ is  $B$-invariant open and dense in $\hat{M}$.
\end{lemma}
\begin{proof}
Clearly $B\cU$ is open and $B$-invariant.
Suppose $\overline{B\cU}\neq\hat{M}$. Let $\cV\neq\emptyset$ be the exterior in $\hat{M}$ of $B\cU$. Then
$\cV$ is also $B$-invariant and open. Since $\cF$ is a minimal hyperbolic lamination it follows that
for each $x\in{M}$ one has that $B\cU\cap{T^1_x\cF}\overset{def}=B\cU_x$ and $\cV\cap{T^1_x\cF}\overset{def}=\cV_x$ are nonempty
disjoint open subsets of the circle fibre $T^1_x\cF$, therefore $C_x\overset{def}=T^1_x\cF-(B\cU_x\cup\cV_x)$ is a closed set in $T^1_x\cF$
with more than one point.  The set $\mathcal{D}\overset{def}=\hat{M}-(B\cU\cup\cV)$ is a nonempty compact set which is $B$-invariant and
it intersects $T^1_x\cF$ in the set $C_x$ with more than one point. Therefore (by Proposition \ref{proposition:more_than_one_point_implies_minimality}) $\mathcal{D}=\hat{M}$. This is a contradiction, hence
$\cV=\emptyset$. Therefore $\cU$ is dense.
\end{proof}

\noi{\bf Proof of Theorem \ref{theorem:transitivity_of_B}:} The proof of  Theorem \ref{theorem:transitivity_of_B} is by the standard Baire category argument. Let $\{{\cU}_i\}_{i\in\N}$ be a countable open base of the topology of $\hat{M}$. Let $\mathfrak{U}_i$ be the set of points  $y\in\hat{M}$ such that the $B$-orbit of $y$ intersects $\cU_i$ ($i\in\N$) then, by lemma 1,
$\mathfrak{U}_i$ is open and dense. By Baire's category theorem the set $\mathfrak{U}\overset{def}=\cap_{i\in\N}{\mathfrak{U}_i}$
is $B$-invariant and residual in $\hat{M}$. The set $\mathfrak{U}$ is dense and consists of dense $B$-orbits. \\

\bigbreak

In the case where the $B$-action on $\hat M$ is minimal,
we have the following:

\begin{proposition}
\label{remark:geodesic_flow_is_topologically_transitive}
Let $(M,\cF)$ be a compact minimal lamination by hyperbolic surfaces. If the $B$-action on $\hat M$ is minimal,
then the laminated geodesic flow on $\hat M$ is topologically transitive.
\end{proposition}

\begin{proof}

Let $U\subset \hat M$ be an open set. It contains a small segment $I$ of a stable horocycle orbit.
 When the $B$-action is minimal, $\cup_{t<0}g_t(I)$ is a dense subset of $\hat M$. One way to see this
 is the following: consider the normalized Lebesgue measure $m$ of the horocycle segment $I$, the sequence of
 probability measures $m_n$ defined by
$\int f\,
 dm_n=\frac{1}{n}\int_0^n \int_{\hat M}f\circ g_t\, dm dt$, and
an accumulation point $m_\infty$ of $m_n$.
In way of contradiction, assume that $K=\overline{\cup_{t<0}g_t(I)}$
is not the whole $\hat M$.
Consider a continuous function
$f:\hat  M\to[0,1]$ taking value 1 on $K$ and value 0 on some point
outside $K$.
 Clearly $\int_{\hat M}fdm_n=1$ for all $n$, and therefore
$\int_{\hat M}fdm_\infty=1$
 But direct computation shows that $m_\infty$ is $B$-invariant.
Thus $m_\infty$ has full support on $\hat M$ and we must have
$\int_{\hat M}fdm_\infty<1$. A contradiction shows that
$\cup_{t<0}g_t(I)$ and hence $\cup_{t<0}g_t(U)$ is dense in $\hat M$.

 Thus the set of points whose geodesic orbit intersects $U$, namely the set
 $\cup_{t\in\R}g_t(U)$, is open and dense. Therefore, there is a residual set of points whose
 geodesic orbit intersects
 every element of a countable basis for the topology of $\hat M$.
\end{proof}

\medskip

The possibility that the $B$-minimal set $\mathcal{M}$ intersects some fibers $T^1_x\cF$ in one point and others in more than one point
cannot be ruled out. This will be shown in the example at the end of this section. Remark that the cardinality of
$T^1_x\cF \cap \mathcal{M}$
remains constant on leaves. Furthermore, since the function
$$x\mapsto \mbox{diameter}(T^1_x\cF \cap \mathcal{M})$$
is upper semicontinuous, the set $\mathcal{R}=\{x\in M/\ \sharp(T^1_x\cF \cap \mathcal{M})=1\}$ is
a countable intersection of open and dense subsets of $M$, unless it is empty.

In order to present the aforementioned example, we will make a small digression to recall some facts about hyperbolic three-manifolds.

Let $\Gamma\subset PSL(2,\R)$ be a cocompact surface group, and $\Sigma=\Gamma\backslash\mathbb{H}$ the hyperbolic surface it defines.
In $\Sigma$, any simple closed curve $\alpha$, not homotopic to a constant, is freely
homotopic to a unique geodesic $C_\alpha$.

Thus, to a free homotopy class of loops
$[\alpha]$ we can associate the hyperbolic length
$\ell([\alpha])=\ell(C_\alpha)$ of
its geodesic representative.

 Let $f:\Sigma\to\Sigma$ be a diffeomorphism
and $f_*$ the induced map on free homotopy classes.
Then, $f$ is said to be of {\it pseudo-Anosov} type if there exist
constants $A([\alpha]), B>0$  such that
$\ell({f_*^n}[\alpha])\geq{A([\alpha]){e^{Bn\ell([\alpha])}}}$,
where $f^n=f\circ\cdots\circ{f}$ ($n$ times).
Thus pseudo-Anosov diffeomorphisms increase exponentially the length of
loops not homotopic to a constant loop.
The property of a diffeomorphism being
of pseudo-Anosov type depends only on the
isotopy class of the diffeomorphism.

Let $f:\Sigma\to\Sigma$ be an orientation-preserving diffeomorphism
of pseudo-Anosov type. Then, by a celebrated
theorem by Thurston (\cite {Gr}, \cite {McM},
\cite {Ot}, \cite{Th1}, \cite{Th2}, \cite{Th3})
the  mapping torus of $f$ is a hyperbolic 3-manifold. Its
hyperbolic metric is unique by Mostow's rigidity theorem. Namely, there is a Klenian
group $\Lambda\subset PSL(2,\mathbb{C})$ such that $\Lambda\backslash\mathbb{H}^3$ is the suspension of $f$.

The suspension construction implies that $\Lambda/\mathbb{H}^3$ fibers over
 $\mathbb{S}^1$ with fiber $\Sigma$. The group $\Lambda$ is a semi-direct product of $\Gamma$ and $\mathbb{Z}$, that is,
 there is an exact sequence
 $$1\to\Gamma\to\Lambda\to\mathbb{Z}\to 0.$$
 Therefore, there is an injective homomorphism $\phi:\Gamma\to\Lambda$ whose image $\phi(\Gamma)$ is a normal subgroup of $\Lambda$.

 Denote by $\mathbb{S}^1$ the circle at infinity of $\mathbb{H}^2$, and by $\mathbb{P}^1$ the sphere at infinity of $\mathbb{H}^3$.
 According to a theorem by Cannon and Thurston (see \cite{Cannon-Thurston}) there is a $\phi$-equivariant continuous and surjective map
 $g:\mathbb{S}^1\to \mathbb{P}^1$ (a sphere filling curve). Let $\Gamma$ act on $\mathbb{S}^1\times \mathbb{P}^1$ by
 $\gamma\cdot(\theta,\zeta)=(\gamma\theta,\phi(\gamma)\zeta)$. Then the graph $\mathcal{X}$ of $g$ is a minimal set of this action.

\begin{example}
\label{example:matsumoto's_counterexample}
\end{example}

 Consider the manifold
 $$M=\Gamma\backslash(\mathbb{H}^2\times\mathbb{P}^1),$$
 where the action is given by $\gamma\cdot(x,\zeta)=(\gamma x, \phi(\gamma)\zeta)$, and the horizontal foliation $\mathcal{F}$ on $M$.
 The leafwise unit tangent bundle of $\cF$ is
 $$\hat M=\Gamma\backslash (PSL(2,\mathbb{R})\times\mathbb{P}^1).$$

 Let $B$ be the affine group. The $B$-action on $\hat M$ is dual (i.e. Morita equivalent) to the $\Gamma$-action on $(PSL(2,\mathbb{R})/B)\times \mathbb{P}^1$ or $\mathbb{S}^1\times\mathbb{P}^1$. There is a minimal set $\mathcal{M}$ of the former $B$-action which corresponds to the above minimal set $\mathcal{X}$ of the latter $\Gamma$-action.

 Let $x\in M$ and let $L$ be the leaf of $\cF$ through $x$.
 Take $\zeta\in\mathbb{P}^1$ such that $\mathbb{H}^2\times\{\zeta\}$ projects onto $L$.
 The intersection $T^1_x\cF \cap \mathcal{M}$ consists of a single point if and only if $(\Gamma\mathcal{M})\cap (PSL(2,\mathbb{R})\times\{\zeta\})$ consists of a single $B$-orbit.
 Since $\Gamma\mathcal{M}=\mathcal{X}B$, this is equivalent to $\mathcal{X}$ intersecting $\mathbb{S}^1\times\{\zeta\}$
 in a single point.  This happens if and only if $\zeta$
has only one preimage under $g$, which is true for some, but certainly
not for any $\zeta$. The last statement is clear since
the map $g$ cannot be a homeomorphism. To show the first statement, let
$\zeta,\zeta'$ be the fixed points of
 some $\phi(\gamma)\in\phi(\Gamma)\setminus\{e\}$, a loxodromic
 element.
 Both $g^{-1}(\zeta)$ and $g^{-1}(\zeta')$,
being closed, nonempty and invariant by $\gamma$, must
contain at least one point from ${\rm Fix}(\gamma)$.
The set $g^{-1}(\zeta)$ cannot contain any point from
$\mathbb S^1\setminus{\rm Fix}(\gamma)$, since if it does, $g^{-1}(\zeta)$
would
contain
the whole ${\rm Fix}(\gamma)$, contradicting the fact that
$g^{-1}(\zeta)\cap g^{-1}(\zeta')=\emptyset$. The same is true
for $g^{-1}(\zeta')$ and
one can conclude that
$g^{-1}(\zeta)$, as well as $g^{-1}(\zeta')$, is a singleton. 

\bigbreak

Theorem \ref{theorem:uniqueness_of_B_minimal_set} implies that the minimal set $\cM$ in Example \ref{example:matsumoto's_counterexample}
is unique.

\bigbreak

We will finally give a condition that guarantees the minimality of the $B$-action on $\hat M$.

\begin{theorem}
Let $(M,\cF)$ be a compact minimal lamination by hyperbolic surfaces that has a holonomy-invariant measure. Then the $B$-action
on its unit tangent bundle is minimal.
\end{theorem}

\begin{proof}
Assume for contradiction that there is a proper minimal set $\mathcal{M}$ for the $B$-action. Since $\cF$ is minimal, we have
$\pi(\mathcal{M})=M$. Consider a foliated chart $D\times T$ of $\cF$, where $T$ is a transversal. Its inverse image can be written
as
$$\Pi^{-1}(D\times T)=D\times \mathbb{S}^1\times T,$$
where $D\times\mathbb{S}^1\times\{t\}$ is a plaque of $T^1\cF$ and $D\times\{\theta\}\times\{t\}$ is a plaque of the foliation by
$B$-orbits. The intersection $\mathcal{M}\cap (D\times\mathbb{S}^1\times T)$ has the form
$$\mathcal{M}\cap (D\times\mathbb{S}^1\times T)=D\times\mathcal{N}$$
for some closed subset $\mathcal{N}$ of $\mathbb{S}^1\times T$. Since $\mathcal{N}$ and $T$ are standard Borel spaces and the
inverse image of any point in $T$ under the projection $\mathcal{N}\to T$ is compact, there is a Borel cross section $\sigma:T\to \mathcal{N}$. This can also be shown in an elementary way as follows. For each $t\in T$, the set
$$\mathcal{N}_t=\{\theta\in \mathcal{S}^1/\ (\theta,t)\in\mathcal{N}\}$$
is a proper nonempty closed subset of $\mathbb{S}^1$. In fact, if it is $\mathbb{S}^1$ for some $t$, $\mathcal{M}$ would contain at least
one leaf of $T^1\cF$. Since $\cF$ is minimal, this would contradict the properness of $\mathcal{M}$. Choosing $T$ smaller if necessary,
one may assume that there is an open interval $I$ of $\mathbb{S}^1$ such that $\mathcal{N}_t$ is contained in $I$. Then there is an upper semi-continuous cross section $\sigma:T\to \mathbb{S}^1$ defined by $\sigma(t)=\sup \mathcal{N}_t$.

Let $\mu$ be an ergodic transverse holonomy-invariant measure. Together with the leafwise hyperbolic volume, it defines an ergodic completely invariant harmonic measure $\lambda$. Proposition \ref{proposition:generic_property_for_harmonic_measures} says that
$\mathcal{N}_t$ must be a singleton for $\mu$-almost all $t$. This gives us the following lemma:

\begin{lemma}
The measure $\sigma_*\mu$ is independent of the Borel cross section $\sigma$.
\end{lemma}

Now the family formed by $\sigma_*\mu$ for each foliaton chart yields a transverse invariant measure of the lamination by $B$-orbits on $\mathcal{M}$, and hence a completely invariant measure on $\mathcal{M}$. But the geodesic flow on $\mathcal{M}$ preserves the transverse measure on one hand, and contracts the leafwise hyperbolic measure on the other hand\footnote{In fact one can similarly show that if a lamination comes from a locally free action of any non-unimodular group, it cannot have holonomy-invariant measures}. A contradiction shows the theorem.

\end{proof}

\section{Minimal sets for the laminated horocycle flow}

We will begin this section looking at several examples.

\begin{example}
\label{example:suspension_of_Anosov_automorphism_of_the_torus}
\end{example}
Let $A:\mathbb{T}^2\to\mathbb{T}^2$ be a hyperbolic linear automorphism of the 2-torus, for example, the one given by the
matrix
$\left(\begin{array}{cc}2&1\\1&1\end{array}\right)$. The suspension of $A$ is a 3-manifold $\mathbb{T}^3_A$ which fibers
over the
circle with fiber the torus $\mathbb{T}^2$. It is a solvmanifold whose universal cover is a solvable Lie group whose Lie
algebra is generated by $X$, $Y$ and $Z$ which satisfy
$$[X,Y]=0,\ [Z,X]=-X, \ [Z,Y]=Y.$$
Therefore, $Z$ and $Y$ generate a locally free action of the affine group, as well as $Z$ and $X$. This is the classical
example of an Anosov flow.
In the unit tangent bundle of the lamination, the $B$-action is not minimal. 
The periodic orbits of the suspension flow are dense in $M=\mathbb T^3_A$.
Each orbit admits two lifts as periodic orbits of the laminated geodesic
flow on $\hat M$, one with the orientation given by the time
parametrization of the suspension flow, and the other opposite.
The former is contained in the unique minimal set $\mathcal M$ of the $B$-action
and, all together, forms a dense subset there. The latter is contained in the
unique minimal set $\mathcal M'$ of the central unstable foliation on $\hat M$.
The set $\mathcal M$ (resp. $\mathcal M'$) is a repellor (an attractor)
of the laminated geodesic flow. Minimal sets for the horocycle flow are 2-tori.

\begin{example}
 \label{example:suspension_of_solenoid}
\end{example}

Let $\mathcal{S}_2$ be the dyadic solenoid
$$
\mathcal{S}_2=\lim_{\longleftarrow}\left\{ \mathbb{S}^1
\overset{f}{\longleftarrow}
\mathbb{S}^1
\overset{f}{\longleftarrow}
\mathbb{S}^1\cdots\right\},$$
where $f(z)=z^2$. Namely,
$$\mathcal{S}_2=\{(z_1,z_2,\ldots):\ z_i\in \mathbb{S}^1, z_{i+1}^2=z_i\ \forall i=1,2,\ldots\}.$$
$\mathcal{S}_2$ is a solenoidal 1-dimensional abelian group, and we consider the map
$$T(s)=s^2,$$
that is, $T(z_1,z_2,\ldots)=(z_1^2,z_1,z_2, \ldots)$, which is an expanding automorphism ($T\circ\theta=\theta\circ T=\mbox{Id}$, where
$\theta$ is the shift map). See \cite{Vietoris}.
The suspension of  $T$ is a lamination $(M,\cF)$, whose 2-dimensional
leaves are the saturation of the leaves of $\mathcal{S}_2\simeq \mathcal{S}_2\times\{0\}$ by the flow. This lamination
is given by orbits of an action of the affine group $B$
for which the orbits of $D$ are the flow lines and the orbits of $U$ are the leaves of the solenoid
$\mathcal{S}_2\times \{t\}$, $t\in [0,1]$. Let $\mathcal{M}\subset \hat M$
be the image of the section $X$ of the bundle $\Pi:\hat M\to M$ defined as follows: for any $x\in M$,
\begin{equation}
\label{equation:section}
X(x)=\left(\frac{d}{dt}\right)_{t=0}\left(x\cdot\left(
\begin{array}{cc}
e^{t/2}&0\\
0&e^{-t/2}
\end{array}
\right)\right)
\end{equation}
Then, as
in the proof of Theorem \ref{proposition:one_point_implies_affine_action}, $\mathcal{M}$
is minimal for the action of $B$ on $\hat M$. It is not minimal for the laminated horocycle flow alone; in fact, each
set of the form $X(\mathcal{S}_2\times \{t\})$, $t\in [0,1]$, is
minimal for the horocycle flow. The restriction to $\mathcal{M}$ of the geodesic flow is the suspension of
$T\circ \pi$ ($\pi$ being, as before, the restriction of the
projection $\Pi$ to $\mathcal{M}$). We observe that the map $T$ has dense periodic orbits.  Hence the suspension flow
has dense periodic orbits. The periodic orbits of the laminated geodesic
flow are just as in Example \ref{example:suspension_of_Anosov_automorphism_of_the_torus}. 

\begin{example}
 \label{example:suspension_of_solenoid_times_minimal_dynamics}
\end{example}

A small twist yields a more general example which is similar in spirit, but without periodic orbits for the geodesic flow.
Let $T:\mathcal{S}_2\to\mathcal{S}_2$ be the expanding map on the dyadic solenoid as in the previous paragraph,
and $F:\mathcal{X}\to\mathcal{X}$ a minimal continuous dynamical system on a compact space $\mathcal{X}$. Now let
$M$ be the suspension of the map
$$
\begin{array}{rcl}
\mathcal{S}_2\times \mathcal{X} &\to& \mathcal{S}_2\times\mathcal{X}\\
(x,y)&\mapsto&(T(x),F(y)).
\end{array}
$$
As before, it carries a minimal 2-dimensional lamination $\cF$ which comes from an action of the affine group.
The image of the section defined by equation \ref{equation:section}
is invariant under the central-stable foliation and is homeomorphic, as a foliated space, to $(M,\cF)$. Minimal
sets for the horocycle flow are sets of the form
$X((\mathcal{S}_2\times\{y\})\times\{t\})$, with $y\in \mathcal{X}$, $t\in [0,1]$. Again $\mathcal{M}$ is the union
of minimal sets for the horocycle flow, and the restriction
of the geodesic flow to $\mathcal{M}$ is a suspension.

In these three examples, the set $\cM$ which is minimal for the action of $B$ is a union of sets which are minimal for
the horocycle flow.
Needless to say, there are examples where the
horocycle flow is itself minimal in $\hat M$, the main one being when the lamination consists of only one leaf which
is a compact
hyperbolic surface. We will now give an example where the horocycle flow is minimal in $\mathcal{M}$ but not in $\hat M$.

\begin{example}
\label{example:central_stable_foliation_of_the_geodesic_flow_of_a_surface}
\end{example}

Let $S$ be an orientable compact connected hyperbolic surface. Then its unit tangent bundle $M=T^1S$ has a locally free
transitive
action of $PSL(2,\R)$ and therefore it has a locally free action of the affine group. Let $\cF$ be the foliation given
by the
orbits of this affine action. Let $(\hat M,T^1\cF)$
be the unit tangent bundle of this foliation. As before, the image of the section defined by
\ref{equation:section} is a minimal set $\mathcal{M}$ for the $B$-action on $\hat M$. This action is obviously
differentiably
conjugate to the original action of the affine group on $M$. Unlike in the previous examples, the set $\mathcal{M}$ is
also minimal for the foliated horocycle flow.

\begin{remark}
When $G=PSL(2,\R)$ acts smoothly and transitively on a 3-manifold $M$ (which must therefore be diffeomorphic to a
quotient of $G$ under
a discrete subgroup), any closed set invariant under the action of the unipotent group $U$ is either a closed orbit or
invariant
under the normalizer of $U$ which is $B$ (see for example
\cite{Witte}). If $M$ is compact, $U$ can have no closed orbits, so any $U$-invariant compact set is also $D$-invariant.
The previous example shows an instance where
this sort of {\it Mautner phenomenon} does not hold when the $G$-action fails to be transitive.
\end{remark}

\bigbreak

In the rest of this section $(M,\cF)$ will be a compact minimal lamination by hyperbolic surfaces such that
the action of $B$ on its unit tangent bundle $(\hat M, T^1\cF)$
{\em is minimal}.

Let $K\subset \hat M$ be a compact invariant set for the laminated horocycle flow $h_t$. It may be that $K=\hat M$.
In any case, since the action of $B$ on $\hat M$ is minimal,
$$\underset{t\in \R}\cup g_t(K),$$
being invariant under both $g_t$ and $h_t$, is dense in $\hat M$.
For every $t\in\R$ the set $g_t(K)$ is also compact and invariant under $h_t$. Assume that
whenever $g_t(K)\cap K\neq\phi$,
in fact $g_t(K)=K$. This holds, if, for example, $K$ is minimal.

Consider the additive subgroup of $\R$ defined as
$$ \mathfrak{G} =\{t\in \R :\ \ g_t(K) = K\}. $$
When the laminated horocycle flow is not minimal in $\hat M$, $\mathfrak{G}$ is either cyclic or trivial. In the first
case,
let $t_0$ be its generator.
The minimality of the affine group action implies that
$$
\underset{t\in\R}\cup g_t(K)=\underset{t\in[0,t_0]}\cup g_t(K)=\hat M,
$$
and $K$ is therefore a global transverse section of the geodesic flow $g_t$, which
every geodesic orbit intersects exactly at intervals of length $t_0$. We call a closed
set having this property a {\em synchronized global transverse section}. The function
$$
\begin{array}{rcl}
 p: \hat M & \to & \mathbb{S}^1=\R/\mathfrak{G}\\
x &\mapsto & t({\rm mod}t_0)
\end{array}
$$
if $g_{-t}(x)\in{K}$ is well defined and it is a locally trivial
fibration of $\hat M$ over $S^1$. That is, the geodesic flow is a
suspension. This was first noticed by Plante in \cite{Plante2}.

The main result of this section states that under the assumption that the $B$-action on $(M, \cF)$ is minimal,
the laminated geodesic flow is never a
suspension.

\begin{proposition}
\label{proposition:g_is_not_a_suspension}
Let $(M,\cF)$ be a compact minimal lamination by hyperbolic surfaces such that the $B$-action on its unit tangent bundle
is minimal.
Then the geodesic flow on $(\hat M, T^1\cF)$ admits no synchronized global transverse section.
\end{proposition}

\begin{proof}
 Suppose $K$ is a synchronized global transverse section for the laminated geodesic flow on $(\hat M, T^1\cF)$.

Let $T:\hat M \to \hat M$ be the involution that leaves every unit tangent space
$T^1_x\cF$, $x\in M$, invariant and takes a unit tangent vector $v$ to $-v$.
We can always assume that $\cF$ is oriented; otherwise we take an orientable
double covering. Under this assumption, $T$ is homotopic to the identity $Id$
in $\hat M$. Let $H_u$, $u\in [0,1]$ be an homotopy taking $T=H_0$ to $Id=H_1$.

There exists an infinite cyclic covering $\psi: \R\times K \to
\hat M$ with the property that $\psi (t+s,x) = g_t \psi(s,x)$ for
all $t,s$. The flow $f_t$ in $\R\times K$ defined by
$f_t(s,x)=(t+s,x)$ is therefore the lifting of $g_t$.

Notice that $T\circ g_t \circ T^{-1} =T\circ g_t \circ T = g_{-t}$, for
all $t\in \R$, and in particular the flows $g_t$ and $g_{-t}$ are
topologically conjugated.

Let us compactify $\R\times K$ by adding two points ``to the left"
and ``to the right". Namely, the compactification is $X=\left
(\R\times K\right ) \cup \{L, R\}$; a neighborhood of $L$ is a set
containing $V_a=\{(t,x):\ t<a\}$, for some $a\in\R$ and
neighborhoods of $R$ are defined analogously.

The flow $f_t$ can be continuously extended to a flow $\bar f_t$ in $X$
that has $L$ and $R$ as fixed points and that satisfies that
\begin{equation}
\label{eq:f_goes_from_left_to_right}
\lim_{t\to -\infty} \bar f_t(x)=L, \ \ \ \lim_{t\to +\infty} \bar f_t(x)=R,
\end{equation}
for every $x\in X\backslash \{L,R\}$.

Likewise, the homotopy $H$ can be lifted to a homotopy $\bar H:~
[0,1]~\times~X~\to~X$ such that $\bar H_1=\bar H (1,\cdot)$ is the
identity in $X$. Then, each map $\bar H_u= \bar H(u,\cdot)$ must
satisfy $\bar H_u(L)=L$, $\bar H_u(R)=R$. Nevertheless, $\bar H_0$
conjugates $f_t$ to $f_{-t}$, which combined with equation
(\ref{eq:f_goes_from_left_to_right}) implies that $H_0(L)=R$ and
$H_0(R)=L$.
\end{proof}

\bigbreak
\begin{remark}
 The proof of the previous proposition in fact proves the following: If a flow is conjugate to its inverse by a
 homeomorphism isotopic to the identity, then it does not admit a global synchronized cross-section.
\end{remark}

\bigbreak

We have the following corollary:

\begin{corollary}
Assume that the $B$-action is minimal on $\hat M$.
If $K\subset \hat M$ is a compact minimal set for the horocycle flow $h_t$, then its intersection
with any central-stable leaf is one of the following:
\begin{itemize}
\item[(i)] the empty set;
\item[(ii)] the whole central-stable leaf;
\item[(iii)] a single horocycle.
\end{itemize}
\end{corollary}

\begin{proof}
 If $o_1$ and $o_2$ are two horocycle orbits on the same central-stable leaf, there is a time $t$ such that
 $g_t(o_1)=o_2$. Let $K$ be a minimal set for the horocycle flow.
 Assume that $K$ intersects a given central-stable leaf $x\!\cdot\! B$ on more than one horocycle orbit.
 If $o_1$ and $o_2$ are two distinct horocycles
 in $(x\!\cdot\! B)\cap K$, there is a $t\neq 0$ such that $g_t(K)\cap K\neq\phi$. Proposition
 \ref{proposition:g_is_not_a_suspension} then says that $K$ must in fact be invariant under the geodesic flow as well,
 and therefore $x\!\cdot\! B\subset K$.
\end{proof}

\bigbreak

As Example \ref{example:suspension_of_solenoid_times_minimal_dynamics} shows, the family of minimal sets for the flow
$h_t$ can in general be very large. Nevertheless,
except for the cases when the $B$-action on $(M,\cF)$ is not minimal, we know of no example of a minimal lamination by
hyperbolic surfaces such that the
horocycle flow $h_t$ is not minimal. Namely, having ruled out the possibility that the group $\mathfrak{G}$ be cyclic,
we do not know if it can ever be trivial.

{\em {\bf Question:} Is it true that for any compact minimal lamination $(M,\cF)$ by hyperbolic surfaces, if
the joint action of the laminated geodesic and horocycle flows is minimal, then the laminated horocycle flow $h_t$
is minimal?}

A positive answer to this question would constitute a generalization of Hedlund's theorem to surface laminations.

For the moment, an interesting corollary of Proposition \ref{proposition:g_is_not_a_suspension} is the following:

\begin{theorem}
If $(M,\cF)$ is a compact minimal lamination by hyperbolic surfaces for which the $B$-action is minimal
and that has a leaf which is not simply connected, then the horocycle
flow $h_t$ is topologically transitive on $\hat M$.
\end{theorem}

\begin{proof}
A leaf of $\cF$ which is not simply connected must have a closed geodesic orbit, since its fundamental group cannot have
elliptic or parabolic elements. Let $T$ be the period
of this closed geodesic orbit and $x$ one of its points. If the orbit of $x$ under the horocycle flow is dense, we have
nothing to prove. Assume it is not dense, and let $K$ be its closure.
Consider the set
$\mathcal{Z}$ of subsets $K'$ of $\hat M$ which are compact and invariant under $h_t$ and
such that $g_T(K')\cap K'\neq \phi$. Since $K\in \mathcal{Z}$, $\mathcal{Z}$ is nonempty.
It follows from Zorn's Lemma that $\mathcal{Z}$ has a minimal element $K_0$, which is contained in $K$. Clearly $g_T(K_0)=K_0$, and $K_0$ is
not invariant under the geodesic flow. It is therefore a global transverse section for $g_t$, which is impossible
according to Proposition \ref{proposition:g_is_not_a_suspension}.
\end{proof}

\begin{remark}
In fact we have proved that if $\cF$ is a compact minimal lamination by hyperbolic surfaces for which the
affine action is minimal, then all periodic points for the geodesic
flow have dense orbits under the horocycle flow.
\end{remark}

\medskip

We will finish this section by showing that the laminated horocycle flow on Sullivan's Universal Hyperbolic Solenoid is
minimal. The Universal Hyperbolic Solenoid is a
compact minimal lamination by hyperbolic surfaces, and its leaves are simply connected.

\begin{example}
 \label{example:Universal_Hyperbolic_Solenoid}
\end{example}

Let $\Sigma_0$ be a compact hyperbolic surface and $x_0$ a point in $\Sigma_0$. We consider the family of all marked
hyperbolic surfaces $(\Sigma,x)$ which are finite regular covers
of $\Sigma_0$ such that the covering map sends $x$ to $x_0$, up to homeomorphisms which send marked points to marked
points. Let
$$\mathcal{C}=\{(\Sigma_\alpha, x_\alpha):\ \alpha \in A\}$$
be this family, and $\pi_\alpha:\Sigma_\alpha\to \Sigma_0$ be the covering map which corresponds to $\alpha \in A$.
A partial order can be defined on $\mathcal{C}$ by
stating that $(\Sigma_\alpha,x_\alpha)\leq (\Sigma_\beta,x_\beta)$ (or $\alpha\leq\beta$ for short) if there exists a
finite regular cover $\pi_{\alpha\beta}:\Sigma_\beta\to\Sigma_\alpha$
such that $\pi_{\alpha\beta}(x_\beta)=x_\alpha$.
The projective limit of $(\mathcal{C},\leq )$ is
$$
\mathcal{H}=\lim_{\longleftarrow} (\Sigma_\alpha,x_\alpha)=\{y=(y_\alpha)\in \prod_{\alpha\in A}\Sigma_\alpha:\
\pi_{\alpha\beta}(y_\beta)=y_\alpha \mbox{ whenever } \alpha \leq\beta\},
$$
seen as a topological subspace of the product $\prod_{\alpha\in A}\Sigma_\alpha$. It is a compact laminated space whose
leaves are dense simply
connected hyperbolic surfaces, see \cite{Sullivan}. It does not depend on the surface $\Sigma_0$, and it is called the
Universal Hyperbolic Solenoid.

For $\alpha, \beta\in A$ such that $\alpha\leq\beta$, we call
$\hat\pi_{\alpha\beta}:T^1\Sigma_\beta\to T^1\Sigma_\alpha$ the map naturally defined by $\pi_{\alpha\beta}$ between the
unit tangent bundles of $\Sigma_\beta$ and $\Sigma_\alpha$,
namely the one given by $\hat\pi_{\alpha\beta}(y,v)=(\pi_{\alpha\beta}(y),d_y\pi_{\alpha\beta}(v))$.

Let $\hat{\mathcal{H}}$ be the unit tangent bundle of $\mathcal{H}$. A point in $\hat{\mathcal{H}}$ is of the form
$z=(z_\alpha)\in \prod_{\alpha\in A}T^1\Sigma_\alpha$
such that
$\hat\pi_{\alpha\beta}(z_\beta)=z_\alpha$ if $\alpha \leq\beta$.

The topology on $\hat{\mathcal{H}}$ has a basis composed of open sets of the form
$$U=\prod_{\alpha\in A}U_\alpha,$$
where $U_\alpha=T^1\Sigma_\alpha$ except for finitely many values of $\alpha\in A$, which we call
$\alpha_1,\ldots,\alpha_n$,
and for each $i$ $U_{\alpha_i}$ is a connected component of $\pi_{\alpha_i}^{-1}(U_0)$, for some fixed small open
set $U_0\subset T^1\Sigma_0$. (The `smallness' of $U_0$ means
that each $\pi_{\alpha_i}$, when restricted to $U_{\alpha_i}$, is a homeomorphism from $U_{\alpha_i}$ to $U_0$.)

We will show that for any $z\in \hat{\mathcal{H}}$ the horocycle through $z$ is dense in $\hat{\mathcal{H}}$, that is,
it intersects every basic open set $U$. Using the notation
introduced in the previous paragraph, let $\beta\in A$ be such that $\alpha_i\leq \beta$ for all $i=1,\ldots n$.
Let $V$ be a connected component of $\pi_\beta^{-1}(U_0)$,
chosen in such a way that $\pi_{\alpha_i\beta}(V)=U_{\alpha_i}$. Since $\Sigma_\beta$ is a compact hyperbolic surface,
Hedlund's theorem tells us that there is a time $t$
for which $h^{(\beta)}_t(z_\beta)\in V$, $h^{(\beta)}_t$ being the horocycle flow on $T^1\Sigma_\beta$. Therefore, at
time $t$ the horocycle orbit of the point $z_{\alpha_i}$
on $T^1\Sigma_{\alpha_i}$ passes through the set $U_{\alpha_i}$, which means that $h_t(z)\in U$.


\bibliography{Hedlund-dec2014}{}
\bibliographystyle{plain}

\end{document}